\documentclass[11pt]{article}

\usepackage{amsmath}
\usepackage{amssymb}
\usepackage{eucal}

\begin{document}

\baselineskip 16pt

\title{Characterizations of some classes of finite  $\sigma$-soluble
 $P\sigma T$-groups\\-{\small Dedicated to  Professor J.C. Beidleman   on the occasion of
 his 80-th birthday}}

\author            
{ Alexander  N. Skiba \\
{\small Department of Mathematics  and Technologies of Programming,}\\  
{\small Francisk Skorina Gomel State University,}\\
{\small Gomel 246019, Belarus}\\
{\small E-mail: alexander.skiba49@gmail.com}}

\date{}
\maketitle

\begin{abstract}      Let   $\sigma =\{\sigma_{i} | i\in I\}$ be some
 partition of the set of all primes $\Bbb{P}$ and  $G$ a  finite group. 
 $G$ is said to be \emph{$\sigma$-soluble} if every chief factor $H/K$ of 
$G$ is a
$\sigma _{i}$-group for some $i=i(H/K)$. 

A set
 ${\cal H}$ of subgroups of $G$  is said to be a  \emph{complete
Hall $\sigma $-set} of $G$  if every member $\ne 1$ of  ${\cal H}$
 is a Hall $\sigma _{i}$-subgroup of $G$ for some $\sigma _{i}\in \sigma $ and
   ${\cal H}$ contains exactly one Hall  $\sigma _{i}$-subgroup of $G$ for every
 $i \in I$ such that $\sigma _{i}\cap   \pi (G)\ne \emptyset$.
A subgroup $A$ of $G$ is said to be   \emph{${\sigma}$-permutable} in $G$
  if $G$ has a complete Hall $\sigma$-set $\cal H$  such that 
 $AH^{x}=H^{x}A$ for all $x\in G$ and all $H\in \cal H$.

We obtain   characterizations of  finite $\sigma$-soluble groups $G$ 
 in which  $\sigma$-permutability    is a transitive relation in $G$.

\end{abstract}


\let\thefootnote\thefootnoteorig

\section{Introduction}

Throughout this paper, all groups are finite and $G$ always denotes
a finite group. Moreover,  $\mathbb{P}$ is the set of all  primes,
  $\pi= \{p_{1}, \ldots , p_{n}\} \subseteq  \Bbb{P}$ and  $\pi' =  \Bbb{P} \setminus \pi$. If
 $n$ is an integer, the symbol $\pi (n)$ denotes
 the set of all primes dividing $n$; as usual,  $\pi (G)=\pi (|G|)$, the set of all
  primes dividing the order of $G$.   
 $G$ is said to be  a $D_{\pi}$-group if $G$ possesses a Hall 
$\pi$-subgroup $E$ and every  $\pi$-subgroup of $G$ is contained in some 
conjugate of $E$.

In what follows, $\sigma$  is some partition of  
$\Bbb{P}$, that is,  $\sigma =\{\sigma_{i} |
 i\in I \}$, where   $\Bbb{P}=\cup_{i\in I} \sigma_{i}$
 and $\sigma_{i}\cap
\sigma_{j}= \emptyset  $ for all $i\ne j$;   $\Pi$ is always supposed
to be a    subset of the set $\sigma$ and $\Pi'= \sigma\setminus \Pi$.

By the analogy with the notation   $\pi (n)$, we write  $\sigma (n)$ to denote 
the set  $\{\sigma_{i} |\sigma_{i}\cap \pi (n)\ne 
 \emptyset  \}$;   $\sigma (G)=\sigma (|G|)$.  $G$ is said to be: 
\emph{$\sigma$-primary} \cite{1} if  $G$ is a $\sigma_{i}$-group for some $i$; 
\emph{$\sigma$-decomposable}  \cite{Shem}  or 
\emph{$\sigma$-nilpotent}   \cite{33} if $G=G_{1}\times \dots \times G_{n}$ 
for some $\sigma$-primary groups $G_{1}, \ldots, G_{n}$; \emph{$\sigma$-soluble} \cite{1}
 if every chief factor of $G$ is $\sigma$-primary;   a \emph{$\sigma$-full group
 of Sylow type}  \cite{1} if every subgroup $E$ of $G$ is a $D_{\sigma _{i}}$-group for every
$\sigma _{i}\in \sigma (E)$.  Note in passing, that every 
$\sigma$-soluble group is a   $\sigma$-full group
 of Sylow type  \cite{2}.

A     set  ${\cal H}$ of subgroups of $G$ is a
 \emph{complete Hall $\sigma $-set} of $G$ \cite{2, commun}  if
 every member $\ne 1$ of  ${\cal H}$ is a Hall $\sigma _{i}$-subgroup of $G$
 for some $\sigma _{i} \in \sigma$ and ${\cal H}$ contains exactly one Hall
 $\sigma _{i}$-subgroup of $G$ for every  $\sigma _{i}\in  \sigma (G)$.

Recall also  that a  subgroup $A$ of $G$ is said to be \emph{${\sigma}$-subnormal}
  in $G$ \cite{1} if   there is a subgroup chain  $$A=A_{0} \leq A_{1} \leq \cdots \leq
A_{n}=G$$  such that  either $A_{i-1} \trianglelefteq A_{i}$ or 
$A_{i}/(A_{i-1})_{A_{i}}$ is  ${\sigma}$-primary 
  for all $i=1, \ldots , n$.

{\bf Definition 1.1.}  We say that  a  subgroup $A$ of $G$ is said to be 
 \emph{$\sigma$-quasinormal} or  
   \emph{$\sigma$-permutable}   in $G$  if $G$ possesses
 a complete Hall $\sigma$-set and $A$ permutes with each   Hall $\sigma _{i}$-subgroup $H$ of
 $G$, that is,  $AH=HA$ for all $i \in I$.

{\bf Remark 1.2.}   Using Theorem B in \cite{1}, it is not difficult to show that 
 if $G$ possesses
 a complete Hall $\sigma$-set  ${\cal H}$ such that $AH^{x}=H^{x}A$ for all
  $H\in {\cal H}$ and  all $x\in G$, then  $A$ is $\sigma$-permutable   in $G$.

{\bf Remark 1.3.} (i) In the classical case
 when $\sigma =\sigma ^{0}=\{\{2\}, \{3\}, \ldots 
\}$: $G$  is $\sigma ^{0}$-soluble (respectively $\sigma ^{0}$-nilpotent) if
 and only if $G$ is soluble 
(respectively nilpotent); 
${\sigma} ^{0}$-permutable subgroups are   also called \emph{$S$-quasinormal} or 
 \emph{$S$-permutable} \cite{prod, GuoII}. A  
 subgroup  $A$ of $G$ is subnormal in $G$ if and only if it is $\sigma 
^{0}$-subnormal in $G$.

(ii) In the other classical case when  $\sigma =\sigma ^{\pi}=\{\pi, \pi'\}$:
 $G$  is $\sigma ^{\pi}$-soluble  (respectively $\sigma ^{\pi}$-nilpotent)
 if and only if $G$ is 
   $\pi$-separable  (respectively $\pi$-decomposable, that is,
 $G=O_{\pi}(G)\times O_{\pi'}(G)$); a subgroup $A$ of a  $\pi$-separable group
 $G$ is $\sigma^{\pi}$-permutable in $G$ if and only if $A$ permutes  with 
all Hall $\pi$-subgroups  and with all Hall $\pi'$-subgroups of $G$.  
A  
 subgroup  $A$ of a  $\pi$-separable group
 $G$ is $\sigma ^{\pi}$-subnormal in $G$ if and only if there is a subgroup
 chain  $$A=A_{0} \leq A_{1} \leq \cdots \leq
A_{n}=G$$  such that   $A_{i}/(A_{i-1})_{A_{i}}$
 is  either  a $\pi$-group or  a $\pi'$-group     for all $i=1, \ldots , n$.

(iii) In fact, in the theory of $\pi$-soluble groups ($\pi= \{p_{1}, \ldots , p_{n}\}$)
 we deal with the  partition 
$\sigma =\sigma ^{0\pi }=\{\{p_{1}\}, \ldots , \{p_{n}\}, \pi'\}$ of $\Bbb{P}$. 
Note that  $G$ is $\sigma ^{0\pi }$-soluble  (respectively $\sigma ^{0\pi }$-nilpotent)
 if and only if $G$ is 
   $\pi$-soluble  (respectively $G=O_{p_{1}}(G)\times \cdots \times
 O_{p_{n}}(G)\times O_{\pi'}(G)$). A subgroup $A$ of a  $\pi$-soluble group
 $G$ is $\sigma^{0\pi }$-permutable in $G$ if and only if $A$ permutes  with 
all Hall $\pi'$-subgroups  and with all Sylow $p$-subgroups  of $G$ for all $p\in \pi$.
Note  also that  a 
 subgroup  $A$ of $G$ is $\sigma ^{0\pi }$-subnormal in $G$ if and only if it is
\emph{$\frak{F}$-subnormal} in $G$ in the sense of Kegel \cite{KegSubn}, where  
$\frak{F}$ is  the class of all $\pi'$-groups. 

We say that  $G$ is a {\sl $P\sigma T$-group} \cite{1} if ${\sigma}$-permutability  
is a transitive relation in $G$, that is, if $K$ is a ${\sigma}$-permutable subgroup
 of $H$ and 
 $H$ is a ${\sigma}$-permutable subgroup of $G$, then  $K$ is a
 ${\sigma}$-permutable subgroup of $G$.   
 In the case when $\sigma =\sigma ^{0}$, a  $P\sigma 
T$-group is  also called a \emph{$PST$-group} \cite{prod}.  Note that if 
$G=(Q_{8}\rtimes C_{3})\curlywedge(C_{7}\rtimes C_{3})$   (see \cite[p. 50]{hupp}), where 
$Q_{8}\rtimes C_{3}=SL(2, 3)$   and  $C_{7}\rtimes C_{3}$ is a non-abelian 
group of order 21, then $G$ is not a $PST$-group but $G$ is  a  $P\sigma T$-group,
 where $\sigma =\{\{2, 3\},  \{2, 3\}'\}$.

The   description of
 $PST$-groups   
 was first obtained by   Agrawal \cite{Agr},
  for the soluble  case, and 
  by Robinson in \cite{217}, for the general case. 
In the   further publications,   authors (see,   for example, the recent papers
 \cite{78}--\cite{rend} and Chapter 2 in \cite{prod}) have found out and  described 
many other   interesting characterizations   of soluble  $PST$-groups.

 The purpose of this 
paper is to study $\sigma$-soluble $P\sigma T$-groups in the most general case
 (i.e., without any restrictions on  $\sigma$). In view of
 Theorem B in \cite{1}, $G$ is a  $P\sigma T$-group if and only if every 
$\sigma$-subnormal subgroup of $G$ is $\sigma$-permutable. 
Being based on this result, here we  prove the following revised  version 
of Theorem A in \cite{1}.

{\bf Theorem A. } {\sl  Let     $D=G^{\frak{N_{\sigma}}}$. If  $G$ is a  $\sigma$-soluble 
 $P\sigma 
T$-group, then  
    the following conditions hold:} 

(i) {\sl $G=D\rtimes M$, where $D$   is an abelian  Hall
 subgroup of $G$ of odd order, $M$ is $\sigma$-nilpotent  and  every element of $G$ induces a
 power automorphism in $D$;  }

(ii) {\sl  $O_{\sigma _{i}}(D)$ has 
a normal complement in a Hall $\sigma _{i}$-subgroup of $G$ for all $i$.}

{\sl Conversely, if  Conditions (i) and (ii) hold for  some subgroups $D$ and $M$ of $G$, then $G$ is  a $P\sigma 
T$-group.}

In this theorem,     $G^{\frak{N_{\sigma}}}$  denotes the \emph{$\sigma$-nilpotent
 residual} of $G$,
 that is,  the intersection of all normal subgroups $N$ of $G$ with 
$\sigma$-nilpotent quotient $G/N$.

{\bf Corollary 1.4.}  {\sl If $G$ is a $\sigma$-soluble $P\sigma 
T$-group, then every quotient and every subgroup of $G$ are $P\sigma 
T$-groups. }

In the 
 case when $\sigma =\sigma ^{0}$,   we get from Theorem 
A the following

{\bf Corollary 1.5} (Agrawal  \cite[Theorem 2.3]{Agr}). {\sl Let 
$D=G^{\frak{N}}$  be 
 the nilpotent residual of $G$.  If  $G$ is a soluble $PST$-group, then  $D$  is  an
abelian  Hall subgroup of $G$ of odd  order and every element
 of $G$ induces a power automorphism 
in  $D$.  }

In the case when $\sigma =\sigma ^{\pi}$, we get from Theorem 
A the following 

{\bf Corollary 1.6.}   {\sl  $G$ is  a    $\pi$-separable  $P\sigma ^{\pi}
T$-group  if and only  if  the following conditions hold:} 

(i) {\sl $G=D\rtimes M$, where $D$ is an abelian  Hall subgroup of $G$ of odd order, 
 $M$ is $\pi$-decomposable 
 and  every element of $G$ induces a
 power automorphism in $D$;  }

(ii) {\sl  $O_{\pi}(D)$ has 
a normal complement in a Hall $\pi$-subgroup of $G$;}

(iii) {\sl  $O_{\pi'}(D)$ has 
a normal complement in a Hall $\pi'$-subgroup of $G$.}

In the case when $\sigma =\sigma ^{0\pi }$, we get from Theorem 
A the following 

{\bf Corollary 1.7.}   {\sl   $G$ is a    $\pi$-soluble  $P\sigma ^{0\pi}
T$-group  if and only  if  the following conditions hold:} 

(i) {\sl $G=D\rtimes M$, where $D$ is an abelian  Hall subgroup of $G$ of odd order, 
 $M=O_{p_{1}}(M)\times \cdots \times O_{p_{n}}(M)\times O_{\pi'}(M)$  
 and  every element of $G$ induces a
 power automorphism in $D$;  }

(ii) {\sl  $O_{\pi'}(D)$ has 
a normal complement in a Hall $\pi'$-subgroup of $G$.}

A natural number $n$ is said to be a \emph{$\Pi$-number}
 if  $\sigma (n)\subseteq \Pi$.  A  subgroup $A$ of $G$ is said to be:
   a \emph{Hall $\Pi$-subgroup} of $G$ \cite{1} if    $|A|$ is  a $\Pi$-number 
 and $|G:A|$ is a $\Pi'$-number; a \emph{$\sigma$-Hall subgroup} of $G$ if $A$
 is a Hall $\Pi$-subgroup of $G$ for some  $\Pi\subseteq \sigma$.

The proof of Theorem A is based on many results and observations of the paper \cite{1}.
We  use also the remarkable result of the paper by Alejandre, Ballester-Bolinches and
 Pedraza-Aguilera
 \cite{al} (see also   Theorem 2.1.8 in \cite{prod})
 that in a soluble  $PST$-group $G$ 
any two isomorphic chief factors are $G$-isomorphic.
Finally, in the proof of Theorem A, the following fact is  useful, which is possibly
 independently interesting.

{\bf Theorem B.} {\sl   Let  $G$ have a normal $\sigma$-Hall  subgroup $D$ 
such that}: (i) {\sl $G/D$ is  a   $P\sigma T$-group, and}  (ii) {\sl every  
$\sigma$-subnormal subgroup of $D$ is normal in $G$. 
  If $G$ is a $\sigma$-full group of Sylow type, then $G$ is a   
$P\sigma T$-group. }

{\bf Corollary 1.8} (See Theorem A in \cite{1}). {\sl  
 Let  $G$ have a normal $\sigma$-Hall  subgroup $D$ such that:}
 (i) {\sl $G/D$ is     $\sigma$-nilpotent, and   }(ii) {\sl   every  
 subgroup of $D$ is normal in $G$.   Then $G$ is a   
$P\sigma T$-group. }

In the case when $\sigma =\sigma ^{0}$, we get from Theorem 
B the following 

{\bf Corollary 1.9} (Agrawal  \cite[Theorem 2.4]{Agr}). {\sl 
  Let  $G$ have a normal Hall  subgroup $D$ such that:} (i) {\sl $G/D$ is a   
$PST$-group, and}  (ii) {\sl every subnormal  
 subgroup of $D$ is normal in $G$.  Then $G$ is a   
$PST$-group. }

Some other applications of Theorems A and B and  some other characterizations
 of $\sigma$-soluble 
$P\sigma T$-groups we discuss in Section 4.

\section{Some preliminary  results}

In view of Theorem  B in \cite{2}, the following fact is true.  

{\bf Lemma 2.1.} {\sl If $G$ is $\sigma$-soluble, then $G$ is a $\sigma$-full group
 of Sylow type.    
}

{\bf Lemma 2.2 } (See  Corollary 2.4 and Lemma 2.5  in \cite{1}).  {\sl  
The class   of all  $\sigma$-nilpotent groups
 ${\mathfrak{N}}_{\sigma}$               is closed under taking  
products of normal subgroups, homomorphic images and  subgroups. Moreover, if  $E$ is a normal 
subgroup of $G$ and  $E/E\cap \Phi (G)$ is $\sigma$-nilpotent, then 
$E$ is $\sigma$-nilpotent.    }

{\bf Lemma 2.3} (See Proposition 2.2.8  in \cite{15}).   {\sl If 
$N$ is a normal subgroup of $G$, then
 $(G/N)^{{\frak{N}}_{\sigma}}=G^{{\frak{N}}_{\sigma}}N/N.$  }

{\bf Lemma 2.4} (See  Knyagina and  Monakhov \cite{knyag}). {\sl
Let $H$, $K$  and $N$ be pairwise permutable
subgroups of $G$ and  $H$ be  a Hall subgroup of $G$. Then $N\cap HK=(N\cap H)(N\cap K).$}

{\bf Lemma 2.5} (See Lemma 2.8 in \cite{1}). {\sl Let  $A$,  $K$ and $N$ be subgroups of  $G$.
 Suppose that   $A$
is $\sigma$-subnormal in $G$ and $N$ is normal in $G$.  }

(1) {\sl If $N\leq K$ and $K/N$ is $\sigma$-subnormal in $G/N$, then $K$
is $\sigma$-subnormal in $G$}.

(2) {\sl $A\cap K$    is  $\sigma$-subnormal in   $K$}.

(3) {\sl If $A$ is a \emph{$\sigma$-Hall subgroup} of $G$, then $A$ is normal in $G$.}

(4) {\sl If $H\ne 1 $ is a Hall $\Pi$-subgroup of $G$  and $A$ is not  a
 $\Pi'$-group, then $A\cap H\ne 1$ is
 a Hall $\Pi$-subgroup of $A$. }

(5) {\sl $AN/N$ is
$\sigma$-subnormal in $G/N$. }

(6) {\sl If  $K$ is a $\sigma$-subnormal subgroup of  $A$,
then $K$ is $\sigma$-subnormal in $G$}.

{\bf Lemma 2.6} (See Lemmas 2.8, 3.1 and Theorem B in 
\cite{1}).   {\sl   Let  $H$,
  $K$ and $R$ be subgroups of  $G$. Suppose that  $H$ is
 $\sigma$-permutable  in $G$ and $R$ is normal in $G$. Then:}

(1) {\sl $H$ is $\sigma$-subnormal in $G$. }

(2)  {\sl The subgroup  $HR/R$ is $\sigma$-permutable  in $G/R$.}

(3) {\sl If $K$ is a $\sigma _{i}$-group, then $K$ is $\sigma$-permutable  in $G$ if and only
 if
 $O^{\sigma _{i}}(G)  \leq N_{G}(K)$}.

(4) {\sl If  $G$ is  
a $\sigma$-full group of Sylow type and  $H\leq K$, then $H$ is  $\sigma$-permutable in $K$.}

(5) {\sl If  $G$ is  
a $\sigma$-full group of Sylow type, $R\leq K$ and $K/R$ is $\sigma$-permutable  in $G/R$,
 then $K$ is
  $\sigma$-permutable  in $G$.}

(6)  {\sl  $H/H_{G}$ is $\sigma$-nilpotent}.

{\bf Lemma 2.7.}   {\sl  The following statements hold:}

(i) {\sl $G$ is a 
 $P\sigma T$-group if and only if every  $\sigma$-subnormal subgroup of 
$G$ is $\sigma$-permutable in $G$. }

(ii) {\sl 
If  $G$ is a 
 $P\sigma T$-group, then every   quotient $G/N$ of $G$ is also a   
$P\sigma T$-group. }

 {\bf Proof.} (i) This follows from Lemmas 2.5(6) and 2.6(1). 

(ii) Let $H/N$ be a $\sigma$-subnormal subgroup of $G/N$. Then 
$H$   is a $\sigma$-subnormal subgroup of $G$ by Lemma 2.5(1), so $H$ is
 $\sigma$-permutable in $G$ by hypothesis and Part (i). Hence $H/N$  is
  $\sigma$-permutable in 
$G/N$ by Lemma 2.6(2). Hence   $G/N$  is  a   
$P\sigma T$-group by Part (i).

The lemma is proved.

\section{Proofs of Theorems A and B}

{\bf Proof of Theorem B.}  Since  $G$ is a $\sigma$-full group of Sylow type by hypothesis,
 it  possesses a 
complete Hall $\sigma $-set   ${\cal H}=\{H_{1}, \ldots , H_{t} \}$ and a
 subgroup $H$ of $G$ is 
$\sigma$-permutable in $G$ if and only if $HH_{i}^{x}=H_{i}^{x}H$ for all
 $H_{i} \in {\cal H}$ and $x\in G$.  
  We can assume without loss  of generality  that $H_{i}$ is a 
$\sigma _{i}$-group for all $i=1, \ldots , t$.

  Assume that
this theorem  is false and let $G$ be a counterexample of minimal order.  
Then  $D\ne 1$ and   for some  $\sigma$-subnormal subgroup $H$ of $G$ and for
some $x\in G$ and $k\in I$
 we have  $HH_{k}^{x}\ne H_{k}^{x} H$ by Lemma 2.7(i).  Let $E=H_{k}^{x}$.

(1) {\sl The hypothesis holds for every quotient  $G/N$ of $G$. }

 It
 is clear that  $G/N$  is a $\sigma$-full group of Sylow type and 
  $DN/N$   is  
a normal $\sigma$-Hall  subgroup of $G/N$.  On the other hand,  
$$(G/N)/(DN/N)  \simeq G/DN \simeq (G/D)/(DN/D),$$ so $(G/N)/(DN/N)$ is a  
$P\sigma T$-group  by Lemma 2.7(ii).  Finally, let $H/N$ be a  $\sigma$-subnormal 
subgroup of $DN/N$.     
Then  $H=N(H\cap D)$  and, by Lemma 2.5(1),  $H$ is $\sigma$-subnormal in $G$. 
  Hence $H\cap D$ is 
$\sigma$-subnormal in $D$ by Lemma 2.5(2), so $H\cap D$ is 
normal in $G$ by   hypothesis. Thus $H/N=N(H\cap D)/N$  is normal in 
$G/N$.   Therefore the  hypothesis holds for  $G/N$.

(2) {\sl $H_{G}=1$. }

Assume that $H_{G}\ne 1$. Clearly,   $H/H_{G}$ is  $\sigma$-subnormal in 
$G/H_{G}$. 
Claim (1) implies that  the   hypothesis holds  for 
$G/H_{G}$, so the choice of $G$ implies that  $G/H_{G}$ is a $P\sigma T$-group. Hence 
$$(H/H_{G})(EH_{G}/H_{G})=
 (EH_{G}/H_{G})(H/H_{G})$$ by Lemma   2.7(i). Therefore  $HE=EH$,  a contradiction. 
 Hence $H_{G}= 1$.

(3) {\sl $DH=D\times H$. }

By Lemma 2.5(2), $H\cap D$ is $\sigma$-subnormal in  $D$. Hence  
$H\cap D$ is normal in $G$ by hypothesis, which implies that $H\cap D=1$ by Claim  (2).  
Lemma 2.5(2) implies also that   $H$ is $\sigma$-subnormal in $DH$. But
 $H$ is a $\sigma$-Hall subgroup of $DH$  since $D$ is a 
$\sigma$-Hall subgroup  of $G$ and $H\cap D=1$.   Therefore $H$ is normal in $DH$ by Lemma 
2.5(3), so $DH=D\times H$.

 {\sl Final contradiction.} Since $D$ is a $\sigma$-Hall subgroup of $G$, then either $E\leq D$
 or $E\cap D=1$. But the former
 case is impossible by Claim (3) since $HE\ne EH$, so  $E\cap D=1$.  Therefore $E$ is a 
$\Pi'$-subgroup of $G$, where $\Pi =\sigma (D)$.  
 By the 
Schur-Zassenhaus theorem, $D$ has a complement $M$ in $G$.  Then $M$ is a 
Hall  $\Pi'$-subgroup of $G $ and so for some $x\in G$ we have $E\leq 
M^{x}$ since $G$ is a $\sigma$-full group of Sylow type.  On the other hand, $H\cap M^{x}$
 is a Hall  $\Pi'$-subgroup of 
$H$ by Lemma 2.5(4) and hence $H\cap M^{x}=H\leq  M^{x}$.  Lemma 2.5(2) implies that $H$ is  $\sigma$-subnormal in $M^{x}$.  But $M^{x}\simeq 
G/D$ is a $P\sigma T$-group by hypothesis, so $HE=EH$ by Lemma 2.7(i). 
This contradiction completes the proof of the theorem.

{\bf Sketch of the proof of  Theorem  A.}   Since $G$ is 
$\sigma$-soluble by hypothesis, $G$ is a $\sigma$-full group of Sylow type by Lemma 2.1. 
 Let  ${\cal H}=\{H_{1}, \ldots , H_{n} \}$ be  
  a  complete Hall $\sigma$-set of $G$.
  We can assume without loss  of generality  that $H_{i}$ is a 
$\sigma _{i}$-group for all $i=1, \ldots , n$.

 First suppose that $G$ is  a $P\sigma 
T$-group. We show that Conditions (i) and (ii) hold for $G$ in this case.  Assume that
this    is false and let $G$ be a counterexample of minimal order.  
Then $D\ne 1.$

(1) {\sl If $R$ is a non-identity normal subgroup of $G$, then
 Conditions} (i) {\sl and} (ii) {\sl hold for  $G/R$} (Since the hypothesis holds for $G/R$ by Lemma
 2.7(ii),  this follows from   the choice of $G$).

(2) {\sl If $E$ is a proper $\sigma$-subnormal subgroup of $G$, then 
 $E^{\frak{N_{\sigma}}}\leq D$ and Conditions} (i) {\sl and} (ii) 
{\sl hold for  $E$}.

Every  $\sigma$-subnormal   subgroup $H$ of $E$ is $\sigma$-subnormal 
in $G$ by Lemma 2.5(6), so    $H$ is  $\sigma$-permutable 
in $G$ by Lemma 2.7(i).     Thus $H$ is $\sigma$-permutable in $E$ by Lemma 
2.6(4). Therefore $E$ is a  $\sigma$-soluble   $P\sigma 
T$-group by Lemma 2.7(i),  so Conditions (i) and (ii) hold for  $E$ by the choice of $G$.  
Moreover,  since $G/D\in {\frak{N_{\sigma}}}$ 
  and ${\frak{N_{\sigma}}}$ is a hereditary  class by Lemma 2.2, 
$E/E\cap D\simeq ED/D\in {\frak{N_{\sigma}}}$  and so    
 $E^{\frak{N_{\sigma}}}\leq 
E\cap D\leq D$. 

(3) {\sl $D$ is nilpotent.}

 (4) {\sl  $D$ is a Hall subgroup of $G$. Hence $D$ has
 a $\sigma$-nilpotent complement $M$ in $G$.}

 Suppose
that this is false and let $P$ be a  Sylow $p$-subgroup of $D$ such
that $1 < P < G_{p}$, where $G_{p}\in \text{Syl}_{p}(G)$.  We can assume 
without loss of generality that $G_{p}\leq H_{1}$.

($a^{0}$)  {\sl    $D=P$ is  a minimal normal subgroup of $G$. }

Let $R$ be a minimal normal subgroup of $G$ contained in $D$. 
 Since
 $D$ is  nilpotent by Claim (3),   $R$ is a $q$-group    for some prime   
$q$. Moreover, 
$D/R=(G/R)^{\mathfrak{N}_{\sigma}}$  is a Hall subgroup of $G/R$ by
Claim (1) and Lemma 2.3.  Suppose that  $PR/R \ne 1$. Then  $PR/R \in \text{Syl}_{p}(G/R)$. 
If $q\ne p$, then    $P \in \text{Syl}_{p}(G)$. This contradicts the fact 
that $P < G_{p}$.  Hence $q=p$, so $R\leq P$ and therefore $P/R \in 
\text{Syl}_{p}(G/R)$ and we again get  that 
$P \in \text{Syl}_{p}(G)$. This contradiction shows that  $PR/R=1$, which implies that 
  $R=P$  is the unique minimal normal subgroup of $G$ contained in $D$. Since $D$ is nilpotent,
 a $p$-complement $E$ of $D$ is characteristic in 
$D$ and so it is normal in $G$. Hence $E=1$, which implies that $R=D=P$.

($b^{0}$) {\sl $D\nleq \Phi (G)$.    Hence for some maximal subgroup
 $M$ of $G$ we have $G=D\rtimes M$  }  (This follows  from ($a^{0}$) and  Lemma 2.2 since $G$
 is not $\sigma$-nilpotent).

($c^{0}$) {\sl If $G$ has a minimal normal subgroup $L\ne D$, then
 $G_{p}=D\times (L\cap G_{p})$.
  Hence $O_{p'}(G)=1$. }

Indeed, $DL/L\simeq D$ is a Hall 
subgroup of $G/L$ by Claim (1) and lemma 2.3. Hence  $G_{p}L/L=DL/L$, so $G_{p}=D\times (L\cap G_{p})$.
 Thus  $O_{p'}(G)=1$ since $D < G_{p}$ by Claim ($a^{0}$).

($d^{0}$)  {\sl   $V=C_{G}(D)\cap M$ is a  normal subgroup of $G$ and 
 $C_{G}(D)=D\times V \leq H_{1}$.  }

In view of  Claims  ($a^{0}$) and  ($b^{0}$),  $C_{G}(D)=D\times V$, where $V=C_{G}(D)\cap M$ 
is a normal  subgroup of $G$. Moreover,  $V\simeq DV/D$ is $\sigma $-nilpotent by
 Lemma 2.2.  Let $W$ be a $\sigma 
_{1}$-complement of $V$. Then $W$  is characteristic in $V$ and so it is normal 
in $G$.    Therefore we have  ($d^{0}$) by Claim ($c^{0}$).

($e^{0}$)  $G_{p}\ne H_{1}$.

Assume that $G_{p}=H_{1}$.  Let $Z$ be a subgroup of order $p$ in $Z(G_{p})\cap D$.
Then, since $ O^{\sigma _{1}}(G)=O^{p}(G)$, $Z$ is normal in 
$G$ by Lemmas 2.6(3) and 2.7(i). Hence $D=Z < G_{p}$  by Claim ($a^{0}$) and so $D < C_{G}(D)$. 
   Then  $V=C_{G}(D)\cap M\ne 1$ is a normal subgroup of $G$ and 
  $V\leq H_{1}=G_{p}$ by Claim ($d^{0}$). Let $L$ be a   minimal
 normal subgroup of $G$ contained in $V$. Then  $G_{p}=D\times L$ is a normal  
elementary abelian subgroup of $G$ by Claim  ($c^{0}$). 
  Therefore every subgroup of $G_{p}$ is 
normal in $G$ by Lemma 2.6(3).  Hence  $|D|=|L|=p$.
Let $D=\langle a \rangle$,  $L=\langle b \rangle$ and $N=\langle ab \rangle$.  
Then $N\nleq D$, so in view of the $G$-isomorphisms
 $$DN/D\simeq N\simeq NL/L= G_{p}/L=DL/L\simeq D $$  we get that 
$G/C_{G}(D)=G/C_{G}(N)$ is a $p$-group since $G/D$ is $\sigma$-nilpotent by Lemma 
2.2.
But then Claim ($d^{0}$) implies that  $G$ is a $p$-group. This 
contradiction shows that we have ($e^{0}$).

{\sl Final contradiction for (4).} In view of Theorem A in \cite{2}, $G$ has a $\sigma 
_{1}$-complement $E$ such that $EG_{p}=G_{p}E$. 
Let $V=(EG_{p})^{{\frak{N}}_{\sigma}}$.  By 
Claim ($e^{0}$), $EG_{p}\ne G$.     On the other hand, since $   D\leq 
EG_{p}$ by Claim ($a^{0}$),  $EG_{p}$ is $\sigma$-subnormal in $G$ by Lemma 
2.5(1). 
 Therefore   Claim (2) implies that    $V$  is a Hall subgroup of $EG_{p}$  and  $V\leq D$, 
 so  for a Sylow 
$p$-subgroup $V_{p}$ of $V$ we have $|V_{p}|\leq |P| < |G_{p}|$. 
Hence  $V$ is a $p'$-group and so   
 $V\leq C_{G}(D)\leq H_{1}=G_{p}$  by Claim ($d^{0}$). Thus   $V=1$.  
Therefore $EG_{p}=E\times G_{p}$ is $\sigma$-nilpotent and so $E\leq C_{G}(D)\leq     
H_{1}$. Hence $E=1$ and  so $ D =1$, a contradiction.  Thus,   
$D$ is a Hall subgroup of $G$.  Hence $D$ has
 a complement $M$ in $G$ by the Schur-Zassenhaus theorem and $M\simeq G/D$
 is $\sigma$-nilpotent by Lemma 2.2.

(5)  {\sl $H_{i}=O_{\sigma _{i}}(D)\times S$ 
 for each $\sigma _{i} \in \sigma (D) $.}

First assume that $N=O^{\sigma 
_{i}}(D)\ne 1$. Since  $D$ is nilpotent by Claim (3), $N$ is a  $\sigma 
_{i}'$-group.  
 Moreover,  $G/N$ is a $P\sigma T$-group by  Lemma 2.7(ii) and  so     
the choice of $G$ implies that  $$H_{i}\simeq H_{i}N/N=(O_{\sigma 
_{i}}(D/N))\times (V/N)=(O_{\sigma _{i}}(D)N/N)\times (V/N).$$ 
  Since  $D$ is a Hall subgroup of $H_{i}$ by Claim (4), 
$DN/N$ is a Hall subgroup of $H_{i}N/N$ and so $V/N$ is a Hall subgroup
 of $H_{i}N/N$. Hence $V/N $  is characteristic in $H_{i}N/N$.  On the 
other hand, since $D/N=(G/N)^{{\frak{N}}_{\sigma}}$ is $\sigma$-nilpotent  
by Lemma 2.2,  
 $H_{i}N/N$ is normal in $G/N$ and so $V/N$ is normal in $G/N$.
 The subgroup $N$ has a complement $S$ in $V$ by 
the Schur-Zassenhaus theorem. Thus $H_{i}\cap V =H_{i}\cap NS=S(H_{i}\cap N)=S
 $ is normal in $H_{i}$.

Now assume  that  $O^{\sigma _{i}}(D)=1$, 
that is, $D$ is a $\sigma _{i}$-group. Then $H_{i}$ is normal in 
$G$, so  all subgroups of $H_{i}$ are 
$\sigma$-permutable  in $G$ by Lemmas 2.5(6), 2.7(i) and hypothesis. 
Since $D$ is a normal Hall subgroup of $H_{i}$, it has a 
complement $S$ in $H_{i}$. Lemma 2.6(3)
 implies that $D\leq   O^{\sigma _{i}}(G)\leq N_{G}(S)$. Hence 
$H_{i}=D\times S$.

(6) {\sl Every subgroup $H$ of $D$ is normal in $G$. Hence every element of
 $G$ induces a power automorphism in $D$. }

Since $D$ is  nilpotent by Claim (3), it is enough to consider 
the case when $H\leq O_{\sigma _{i}}(D)=H_{i}\cap D$ for some $\sigma _{i}\in \sigma (D)$.
Claim (5) implies that $H_{i}=O_{\sigma _{i}}(D)\times S$. 
It is clear that $H$ is subnormal in $G$, so $H$ is $\sigma$-permutable in 
$G$. Therefore $$G=H_{i}O^{\sigma _{i}}(G)=
(O_{\sigma _{i}}(D)\times S)O^{\sigma _{i}}(G)=SO^{\sigma _{i}}(G)\leq 
N_{G}(H)$$  by Lemma 2.6(3).

 (7) {\sl  If  $p$ is a  prime such that $(p-1, |G|)=1$, then  $p$
does not divide $|D|$. Hence the smallest prime in $\pi (G)$ belongs to $\pi (|G:D|)$.
 In particular,  $|D|$ is odd. }

Assume that this is false.
 Then, by Claim (6),  $D$ has a maximal subgroup $E$ such that
$|D:E|=p$ and  $E$ is normal in $G$. It follows that  $C_{G}(D/E)=G$ since $(p-1, 
|G|)=1$.   
Hence 
$G/E=(D/E)\times (ME/E)$, where
$ME/E\simeq M\simeq G/D$ is $\sigma$-nilpotent. Therefore $G/E$ is
$\sigma$-nilpotent. But then $D\leq E$, a contradiction. 
Hence we have (7).

(8)  {\sl  $D$ is abelian.}

In view of Claim 
(6), $D$ is a Dedekind group.  Hence $D$ is abelian since $|D|$ is  odd  by Claim (7).  

 From Claims (4)--(8) we get that Conditions (i) and (ii) hold  for $G$.

Now we show that if Conditions (i) and (ii) hold for  $G$,
 then
 $G$ is a $P\sigma T$-group.  Assume that
this   is false and let $G$ be a counterexample of minimal order.  
  Then  $D\ne 1$ and, by Lemma 2.7(i), for some  $\sigma$-subnormal subgroup $H$ of $G$ and for
some $x\in G$ and $k\in I$
 we have  $HH_{k}^{x}\ne H_{k}^{x} H$. Let $E=H_{k}^{x}$.

($1^{0}$) {\sl   If  $N$ is a minimal normal subgroup of $G$, then
 $G/N$ is a $P\sigma T$-group} (Since the hypothesis holds for $G/N$, this follows from
 the choice of $G$).

$(2^{0})$ {\sl If $N$ is a minimal normal subgroup of $G$, then $EHN$ is
 a subgroup of $G$.  Hence $E\cap N=1$.}

Claim $(1^{0})$ implies that $G/N$ is a $P\sigma T$-group. On the other hand,
 $EN/N$ is a Hall $\sigma _{k}$-subgroup of $G/N$ and,
 by Lemma 2.5(5), 
    $HN/N$ is a $\sigma$-subnormal subgroup of 
$G/N$.   Note also that $G/N$ is $\sigma$-soluble, so every two Hall 
$\sigma _{k}$-subgroups of $G/N$ are conjugate by Lemma 2.1. Thus,
 $$(HN/N)(EN/N)=(EN/N)(HN/N)=EHN/N$$ by Lemma 2.7(i).  Hence $EHN$ is a subgroup of $G$.    
Since $G$ is $\sigma$-soluble, $N$ is a $\sigma _{j}$-group for some  $j$. 
Hence in the case   $E\cap N\ne 1$ we have $j=k$, so $N\leq E$. But then   
 $EHN=EH=HE$, a contradiction. Thus $E\cap N=1$.

$(3^{0})$ $|\sigma (D)| > 1$.

Indeed, suppose that $\sigma (D)=\{\sigma _{i}\}$. Then  $H_{i}/D$ is normal in $G/D$ since 
$G/D\simeq M$   is $\sigma$-nilpotent by hypothesis, so  $H_{i}=D\times S$ 
 is normal in $G$.  The subgroup 
 $S$ is also 
 normal in $G$ since it
 is characteristic in $H_{i}$. On the other hand,
 Theorem B  and the choice of $G$ 
imply that $S\ne 1$.

 Let $R$ and $N$ be minimal normal subgroups of $G$ such that  $R\leq 
D$ and $N\leq S$. 
 Then $R$ is a group of order $p$ for some prime $p$  
and $N$ is a $p'$-group since $D$ is a Hall subgroup of $H_{i}$.
 Hence $R\cap HN\leq O_{p}(HN)\leq P$, where $P$ is a Sylow 
$p$-subgroup of $H$, so $R\cap HN=R\cap H$.  
 Claim $(2^{0})$   implies that $EHR$ and   $EHN$ are subgroups of $G$.  
 Therefore from Lemma 2.4 and Claim ($2^{0}$)  we get  that
 $R\cap EHN=R\cap E(HN)=(R\cap E)(R\cap HN)=R\cap H$. 
 Hence  $$EHR\cap EHN=E(HR\cap EHN)=EH(R\cap EHN)$$$$=EH(R\cap H)=EH$$  is a subgroup of $G$. Thus   
$HE=EH$, a contradiction. Hence we have $(3^{0})$.

 {\sl Final contradiction.} Since   $|\sigma (D)| > 1$ by Claim ($3^{0}$)
 and $D$
 is nilpotent, $G$  
has at least two    minimal normal subgroups $R$ and $N$ such that $R, N\leq D$ 
and $\sigma (R)\ne \sigma (N)$.  Then  at least one of the subgroups 
$R$ or $N$, $R$ say, is a $\sigma _{i}$-group for some $i\ne k$. Hence  
 $R\cap HN\leq O_{\sigma _{i}}(HN)\leq V$, where $V$ is a Hall $\sigma _{i}$-subgroup of $H$,
 since $N$ is a $\sigma _{i}'$-group and $G$ is a $\sigma$-full group of 
Sylow type.  Hence  $R\cap HN=R\cap H$.
Claim ($2^{0}$) implies that   $EHR$ and $EHN$ are subgroups of $G$. 
Now,  arguing similarly as in the proof of $(3^{0})$, one can show that 
  $EHR\cap EHN=EH=HE$. 
This contradiction  
 completes the proof of the fact that $G$ is a $P\sigma T$-group.

The theorem is proved.

\section{Some other characterizations of $\sigma$-soluble $P\sigma T$-groups}

  Theorem A and Theorem B in \cite{1}  are     basic in the sense that  many other
 characterizations of $\sigma$-soluble
 $P\sigma T$-groups  
can be obtained by using these two results.  As
 a partial illustration to this, we give in this section our next three characterizations
 of $\sigma$-soluble $P\sigma T$-groups.

1. Recall that   $Z_{\sigma }(G)$ denotes the
 \emph{$\sigma$-hypercentre } of $G$ \cite{ProblemII}, that is, the largest normal subgroup of 
$G$ such that for every chief factor $H/K$ of $G$  below $Z_{\sigma 
}(G)$ the semidirect product $(H/K)\rtimes (G/C_{G}(H/K))$  is $\sigma$-primary. 
 
We say, following \cite[p. 20]{prod}, that a subgroup $H$ of $G$ is 
\emph{$\sigma$-hypercentrally embedded in $G$}, if $H/H_{G}\leq Z_{\sigma }(G/H_{G})$.

{\bf Theorem 4.1.}   {\sl Let  $G$ be  $\sigma$-soluble. Then $G$ is a  
 $P\sigma T$-group if and only if
 every $\sigma$-subnormal subgroup of $G$ is $\sigma$-hypercentrally embedded in $G$.}

{\bf Proof. }  Let $D=G^{{\frak{N}}_{\sigma}}$.  
 First we show that if   $G$ is  a   $P\sigma 
T$-group, then every $\sigma$-subnormal subgroup  $H$ of $G$ is
 $\sigma$-hypercentrally embedded in $G$.  Assume that this is false 
and let $G$ be a counterexample with $|G|+|H|$ minimal.  Then $G/H_{G}$ is a 
  $\sigma$-soluble $P\sigma 
T$-group by Lemma 2.7(ii) and $H/H_{G}$  is $\sigma$-subnormal in 
$G/H_{G}$ by Lemma 2.5(5).  Hence the choice of $G$ implies that  
$H_{G}=1$, so 
$H$ is  $\sigma$-nilpotent by Lemma 2.6(6). Therefore  every subgroup of $H$ is 
$\sigma$-subnormal in $G$ by Proposition 2.3 in \cite{1} and  Lemma 2.5(6).
  Assume that $H$ possesses two distinct
 maximal subgroups $V$ and $W$.  Then $V, W\leq Z_{\sigma }(G)$ by 
minimality of $|G|+|H|$ since $V_{G}=1=W_{G}$, which implies that $H\leq Z_{\sigma }(G)$.
Hence $H $ is a cyclic $p$-group for some  $p\in \sigma _{i}$.

  By Theorem A, $G=D\rtimes M$,
 where $D$   is 
a  Hall
 subgroup of $G$, $M$ is $\sigma$-nilpotent  and  every subgroup of $D$
 is normal in $G$.   Then    $H\cap D=1$ and so, in view of Lemma 2.1, we can
 assume without loss of generality  that  $H\leq M$.
 Lemma 2.7(i) 
implies that  $H$ is $\sigma$-permutable in $G$, so 
$$H^{G}=H^{DM}=H^{O^{\sigma _{i}}(G)M}=H^{M}\leq M$$ by Lemma  2.6(3).
Hence   $H^{G}\cap D=1$ and  then, from the $G$-isomorphism 
$H^{G}D/D\simeq H^{G}$,
 we deduce that $H\leq H^{G}\leq Z_{\sigma }(G)$. Therefore $H$ is 
$\sigma$-hypercentrally embedded in $G$.  This contradiction completes the 
proof of the necessity of the condition of the theorem.

{\sl Sufficiency.} It is enough to show that if  a $\sigma$-subnormal subgroup
 $H$ of
a  $\sigma$-soluble group $G$ is 
$\sigma$-hypercentrally embedded in $G$, then $H$ is $\sigma$-permutable in $G$. 
Assume that this is false  and let $G$ be a counterexample with $|G|+|H|$ 
minimal.     
Since $G$ is $\sigma$-soluble, it is a $\sigma$-full group of Sylow type 
by Lemma 2.1.  
 Therefore, in view of Lemma 2.6(5), 
$H_{G}=1$ and  so   $H\leq Z_{\sigma }(G)$. It is clear that $Z_{\sigma }(G)$ 
is $\sigma $-nilpotent, so 
 $H=H_{1}\times \cdots \times H_{t}$ for some   
$\sigma$-primary groups  $H_{1},  \ldots ,  H_{t}$.   Moreover, the 
minimality of $|G|+|H|$ implies that  $H=H_{1}$ is a $\sigma _{i}$-group for 
some $i$. Hence $H\leq N$, where $N$ is a Hall $\sigma _{i}$-subgroup of $Z_{\sigma }(G)$.
Since  $Z_{\sigma }(G)$ is $\sigma $-nilpotent,  $N$ is characteristic in
 $Z_{\sigma }(G)$  and so  $N$ is normal in $G$.

Let $1= Z_{0} < Z_{1}
< \cdots <  Z_{t} =N$ be a chief
 series of $G$  below $N$ and  $C_{i}=
C_{G}(Z_{i}/Z_{i-1})$. Let $C= C_{1} \cap \cdots \cap C_{t}$.  Then $G/C$ 
 is a ${\sigma _{i}}$-group. On the other hand,
 $C/C_{G}(N)\simeq A\leq \text{Aut}(N)$ stabilizes     
the series    $1= Z_{0} < Z_{1}
< \cdots <  Z_{t} = N$, so   $C/C_{G}(N)$ is a $\pi (N)$-group by  \cite[Ch. A, 12.4(a)]{DH}.
   Hence $G/C_{G}(N)$ is  a ${\sigma _{i}}$-group and  
 so  $O^{\sigma _{i}}(G)\leq C_{G}(N)$.   
But then $O^{\sigma _{i}}(G)\leq  C_{G}(H)$, so $H$ is $\sigma$-permutable in $G$ by 
Lemma 2.6(3).   This contradiction completes the 
proof of the sufficiency of the condition of the theorem.

 The theorem is proved.

In the 
 case when $\sigma =\sigma ^{0}$, we have $Z_{\sigma }(G)=Z_{\infty }(G)$.
 Hence from Theorem 4.1   we get 

{\bf Corollary 4.2} (See Theorem 2.4.4 in \cite{prod}).  {\sl  Let
  $G$ be  soluble. Then $G$ is a   $PST$-group
 if and only if 
 every subnormal subgroup $H$ of $G$ is hypercentrally embedded in $G$ (that is,
 $H/H_{G}\leq Z_{\infty }(G/H_{G}))$.}

2. We say, following \cite[p. 68]{prod}, that $G$ \emph{satisfies property
 ${\cal Y}_{\sigma _{i}}$} if whenever $H\leq 
K$ are two $\sigma _{i}$-subgroups of $G$, $H$ is $\sigma$-permutable in $N_{G}(K)$.

The idea of the next theorem goes back to the paper \cite{A}.

{\bf Theorem 4.3.} {\sl Let  $G$ be  $\sigma$-soluble. Then $G$ is a   $P\sigma T$-group
 if and only if
$G$ satisfies ${\cal Y}_{\sigma _{i}}$ for all primes $i$.}

{\bf Lemma 4.4.} {\sl   Let  $K\leq H$  and $N$ be subgroups of  $G$. 
Suppose that  $K$ is 
 $\sigma$-permutable  in $H$ and $N$ is normal in $G$. Then $KN/N$
 is $\sigma$-permutable  in $HN/N$.}

{\bf Proof.}  Let  $f:H/H\cap N \to HN/N$  be the canonical isomorphism from
$H/H\cap N $ onto $HN/N$. Then $f(K(H\cap N)/(H\cap N))=KN/N,$  so   $KN/N$
 is $\sigma$-permutable  in $HN/N$ by Lemma 2.6(2). 

The lemma is proved.

{\bf Sketch of the proof of Theorem 4.3. } {\sl Necessity.}  Let  $H\leq 
K$ be  two $\sigma _{i}$-subgroups of $G$ and $N=N_{G}(K)$. Then $H$ is
 $\sigma$-subnormal in $N$ by Lemma 2.5(6). On the other hand, Corollary 1.2 implies that $N$ 
is a $\sigma$-soluble  $P\sigma T$-group. Therefore  $H$ is 
$\sigma$-permutable in $N$ by Lemma 2.7(i).

{\sl Sufficiency. } It is enough to show that 
 Conditions (i) and  (ii) of Theorem A  hold for  $G$.  Assume that this is false 
and let $G$ be a counterexample of minimal order. Since $G$ is  
$\sigma$-soluble, it is a $\sigma$-full group of Sylow type by Lemma 2.1.  
Let $D=G^{{\frak{N}}_{\sigma}}$.

(1) {\sl Every proper subgroup $E$ of  $G$   is a $\sigma$-soluble     $P\sigma 
T$-group and $E^{\frak{N_{\sigma}}}\leq D$} (This follows from Lemmas  2.2, 2.3, 2.6(4)
 and the choice of $G$).

(2) {\sl  $G/N$  is a    $\sigma$-soluble  $P\sigma 
T$-group for every minimal normal subgroup $N$ of $G$. }

Let $H/N\leq 
K/N$ be  two $\sigma _{i}$-subgroups of $G/N$.  Since  $G$   is  
$\sigma$-soluble, $N$ is a  $\sigma _{j}$-subgroup for some $j$. Assume 
that $j\ne i$. Then there are a Hall $\sigma _{i}$-subgroup $V$ of $H$ and 
a Hall $\sigma _{i}$-subgroup $W$ of $K$ such that $V\leq W$ since $G$ is a $\sigma$-full group of 
Sylow type. Then $V$ is  $\sigma$-permutable in $N_{G}(W)$ by hypothesis, so $H/N=VN/N$  
is $\sigma$-permutable in $N_{G}(W)N/N=N_{G/N}(WN/N)=N_{G/N}(K/N) $ by Lemma 
4.4.   Similarly we get that  $H/N$  is  $\sigma$-permutable in 
$N_{G/N}(K/N) $   in the case when $j=i$. 

(3) {\sl $D$ is nilpotent.}

(4) {\sl  $D$ is a Hall subgroup of $G$ and $H_{i}=O_{\sigma _{i}}(D)\times S$ 
 for each $\sigma _{i} \in \sigma (D) $} (See  Claims (4) and (5)
 in  the proof of Theorem A and use Claims (1), (2) and (3)).

(5) {\sl Every subgroup $H$ of $D$ is normal in $G$. Hence every element of
 $G$ induces a power automorphism in $D$. }

Since $D$ is  nilpotent by Claim (3), it is enough to consider 
the case when $H\leq O_{\sigma _{i}}(D)=H_{i}\cap D$ for some $\sigma _{i}\in \sigma (D)$.
Hence 
 $H$ is $\sigma$-permutable in 
$G$ by hypothesis.
Claim (4) implies that $H_{i}=O_{\sigma _{i}}(D)\times S$.  Therefore $G=H_{i}O^{\sigma _{i}}(G)=SO^{\sigma _{i}}(G)\leq 
N_{G}(H)$  by Lemma 2.6(3). 

(6) {\sl $D$ is abelian of odd order}  (See  Claims (7) and (8)
 in  the proof of Theorem A and use Claim (5)).

The theorem is proved.

{\bf Corollary 4.5} (Ballester-Bolinches and   Esteban-Romero \cite{A}, see also Theorem 2.2.9
 in \cite{prod}).  {\sl  $G$ is
 a soluble $PST$-group if and only if  $G$ satisfies ${\cal Y}_{p}$ for all primes $p$.}

{\bf Proof.}  It is enough to note that, as  it was remarked at the beginning of the proof
 of Theorem 
2.2.9 in \cite{prod},  every group which satisfies ${\cal Y}_{p}$ for all
 primes $p$  is soluble.

3.  We say that a subgroup $A$ of  $G$ is  \emph{$\sigma$-modular}     
(\emph{$S$-modular}  in the case  $\sigma =\sigma ^{0}$) provided   
$G$ possesses a complete Hall $\sigma $-set and  
 $\langle A, H \cap C \rangle=\langle A, H \rangle \cap C$ for every Hall 
$\sigma _{i}$-subgroup $H$ of $G$ and all $i\in I$ and   $A\leq 
 C \leq   G$.

{\bf Theorem  4.6.} {\sl  Let $G$ be  $\sigma$-soluble. Then $G$ is a  $P\sigma T$-group
 if and only if every  $\sigma$-subnormal subgroup $A$ of $G$ is  
$\sigma$-modular in every subgroup of $G$ containing $A$.  }

{\bf Proof.}  Since  $G$ is  $\sigma$-soluble, 
$G$ is a $\sigma$-full group of Sylow type by Lemma 2.1. Hence 
 $G$  possesses a 
complete Hall $\sigma $-set   ${\cal H}=\{H_{1}, \ldots , H_{t} \}$ and, for each $\sigma _{i} 
\in \sigma (G)$,    a
 subgroup $H$ of $G$ is a Hall  $\sigma _{i}$-subgroup of $G$ if and only 
if $H=H_{k}^{x}$ for some $x\in G$ and  $H_{k}\in {\cal H}$.

{\sl Sufficiency.} Assume that this is false and let $G$ be a 
counterexample of minimal order. Then, in view of Lemma 2.7(i), $G$   has a 
 $\sigma$-subnormal subgroup $A$ which  is not  
$\sigma$-permutable in $G$.    Hence,  for some 
 $H_{i} \in {\cal H}$ and $x\in G$, we have $AH_{i}^{x}\ne H_{i}^{x}A$. 
Note also that every proper $\sigma$-subnormal subgroup $E$ of $G$ is a 
$P\sigma T$-group. Indeed, $E$ is clearly $\sigma$-soluble and
 if $H$ is  a $\sigma$-subnormal subgroup of $E$, then 
$H$ is $\sigma$-subnormal in $G$ by Lemma 2.5(6). Hence  $H$ is $\sigma$-modular in
 every subgroup of $E$ containing $H$ by hypothesis. Thus the hypothesis 
holds for $E$ and so $E$  is a $P\sigma T$-group by the choice of $G$.

By definition, there is a subgroup chain  $A=A_{0} \leq A_{1} \leq \cdots \leq
A_{n}=G$  such that  either $A_{i-1}\trianglelefteq A_{i}$ or $A_{i}/(A_{i-1})_{A_{i}}$
 is  ${\sigma}$-primary for all $i=1, \ldots , n$.  We can assume without 
loss of generality that $M=A_{n-1} < G$. Then $M$  is a $P\sigma T$-group since $M$ is clearly 
$\sigma$-subnormal in $G$, so 
$A$ is $\sigma$-permutable 
in $M$ by Lemma 2.7(i).  Moreover, the $\sigma$-modularity of $A$ in $G$ 
implies that $$M=M\cap \langle A, H_{i}^{x}\rangle=\langle A, (M\cap 
H_{i}^{x})\rangle .$$ On the other hand, by Lemma 2.5(4),  $M\cap H_{i}^{x}$ is a Hall $\sigma 
_{i}$-subgroup of $M$, where  $\{\sigma 
_{i}\}=\sigma (H_{i})$. Hence  $M= A  (M\cap H_{i}^{x})=(M\cap H_{i}^{x})A$. 
 If  $H_{i}^{x}\leq M_{G}$, then  
$A  (M\cap H_{i}^{x})=A H_{i}^{x}=H_{i}^{x}A$ and so $H_{i}^{x}\nleq 
M_{G}$. 

 Now note that  $H_{i}^{x}M=MH_{i}^{x}$. Indeed, if $M$ is normal in $G$, 
it is clear. Otherwise, $G/M_{G}$ is $\sigma$-primary and so $G=MH_{i}^{x}=H_{i}^{x}M$ since 
$H_{i}^{x}\nleq 
M_{G}$ and $H_{i} \in {\cal H}$.  Therefore
 $$H_{i}^{x}A=H_{i}^{x}(M\cap H_{i}^{x})A=H_{i}^{x}M=MH_{i}^{x}    
=H_{i}^{x}(M\cap H_{i}^{x})A= H_{i}^{x}A.$$ This contradiction completes the 
proof of the sufficiency of the condition of the theorem. 

 {\sl Necessity.} In view of Lemma 2.6(4), it is enough to show that if  $A$ is 
 a $\sigma$-subnormal subgroup  of $G$, then  $A$  is $\sigma$-modular in 
$G$. First note that  
$A$ is $\sigma$-permutable in $G$ by Lemma 2.7(i).  Therefore for every 
subgroup $C$ of $G$ containing $A$,  for every $i\in I$, and for all
Hall  
$\sigma _{i}$-subgroup $H$ of $G$  we have 
 $$\langle A, H \cap C \rangle=  A(H\cap C)=AH\cap C= 
\langle A, H \rangle \cap C,$$ so $A$ is  $\sigma$-modular in $G$. 

 The 
theorem is proved.

From Theorem 4.6 we get the following  characterization of 
soluble $PST$-groups.

{\bf Corollary  4.7.} {\sl  Let $G$ be soluble.   Then $G$ is a  $PST$-group
 if and only if every  subnormal subgroup $A$ of $G$ is  
$S$-modular in every subgroup of $G$ containing $A$. }

\end{document}